\documentclass{gtart}

\input gtoutput
\volumenumber{3}\papernumber{7}\volumeyear{1999}
\pagenumbers{155}{165}\published{2 July 1999}
\proposed{Haynes Miller}\seconded{Ralph Cohen, Gunnar Carlsson}
\received{17 March 1998}\revised{27 May 1999}
\accepted{30 June 1999}

\usepackage{amscd,amsmath}

\def\xrightarrow#1{\buildrel{#1\,}\over\to}  

\numberwithin{equation}{section}
\newtheorem{theorem}[equation]{Theorem}
\newtheorem{lemma}[equation]{Lemma}
\theoremstyle{definition}
\newtheorem{remark}[equation]{Remark}
\newtheorem{definition}[equation]{Definition}
\newtheorem{condition}[equation]{Condition}
\DeclareMathOperator{\im}{im}
\DeclareMathOperator{\Ext}{Ext}
\newcommand{\Smash}{\wedge}
\newcommand{\dual}[1]{D \mspace{-2mu} {#1}}
\newcommand{\invlim}{\varprojlim}
\newcommand{\HFp}{H\!\mathbf{F}_{p}}
\newcommand{\Z}{\mathbf{Z}}
\newcommand{\cohomology}{\Ext^{**}_{A}(\ztwo,\ztwo)}
\newcommand{\ztwo}{\mathbf{F}_{2}}
\newcommand{\conn}[1]{\lvert #1 \rvert}
\newcommand{\test}{\mathcal{T}}
\newcommand{\exactcouple}[6]
 {\begin{picture}(16,8)(0,1.7)	   
  \put(2,8){\ensuremath{#1}}	   
  \put(8,2){\ensuremath{#2}}	   
  \put(14,8){\ensuremath{#3}}	   
  \put(5.2,4.8){\makebox(0,0)[r]{\ensuremath{#4}}}
  \put(11.8,4.8){\makebox(0,0)[l]{\ensuremath{#5}}}
  \put(8.5,9){\makebox(0,0){\ensuremath{#6}}}
  \put(3.5,7.5){\vector(1,-1){4.5}}
  \put(9,3){\vector(1,1){4.5}}
  \put(13.5,8.2){\vector(-1,0){9.5}}
 \end{picture}}
\newcommand{\ulp}{\textup{(}}
\newcommand{\urp}{\textup{)}}
\newcommand{\uc}{\textup{:}}
\newcommand{\mathcolon}{\colon\,}

\begin{document}

\title
{Vanishing lines in generalized Adams spectral\\sequences are generic}
\shorttitle{Vanishing lines in Adams spectral sequences}

\author{M\kern .17em J Hopkins\\J\kern .17em H Palmieri\\
J\kern .17em H Smith}
\asciiauthors{MJ Hopkins, JH Palmieri, JH Smith}

\address{Department of Mathematics, Massachusetts Institute of Technology\\
Cambridge, MA 02139, USA\\\smallskip\\
Department of Mathematics, University of Notre Dame\\
Notre Dame, IN 46556, USA\\\smallskip\\
Department of Mathematics, Purdue University \\
West Lafayette, IN 47907, USA}
\email{mjh@math.mit.edu, palmieri@member.ams.org, jhs@math.purdue.edu}
\asciiaddress{Department of Mathematics, 
Massachusetts Institute of Technology
Cambridge, MA 02139, USA
Department of Mathematics, University of Notre Dame
Notre Dame, IN 46556, USA
Department of Mathematics, Purdue University 
West Lafayette, IN 47907, USA}

\keywords{Adams spectral sequence, vanishing line, generic}

\primaryclass{55T15}\secondaryclass{55P42}

\begin{abstract}
We show that in a generalized Adams spectral sequence, the presence
of a vanishing line of fixed slope (at some term of the spectral
sequence, with some intercept) is a generic property.
\end{abstract}

\maketitlepage

\section{Introduction}
\label{sec-intro}

Our main result, Theorem~\ref{thm-main}, was motivated by several
questions.

The thick subcategory theorem of Hopkins and Smith \cite[Theorem
7]{hopkins-smith} describes some important structure in stable
homotopy theory: one can get a good deal of information about a finite
$p$--local spectrum $X$ by knowing its ``type''---that is, knowing the
unique number $n$ so that its $n$th Morava $K$--theory $K(n)_{*}X$ is
nonzero while $K(n-1)_{*}X = 0$.  While one cannot in general
determine the type of a finite $p$--local spectrum by knowing just its
mod $p$ cohomology, one might ask whether the type is reflected
somehow in the behavior of the classical Adams spectral sequence.  Our
main theorem shows that the answer is yes, and we explain this
behavior (the existence of a vanishing line with a particular slope)
as our first application after the statement of the theorem.

Spectral sequences are often sequences of bigraded abelian groups
$\{E_{r}^{**}\}$, indexed starting at $r=1$ or $r=2$.  The
$E_{r+1}$--term is a subquotient of the $E_{r}$--term, so if one has
some sort of vanishing result at $E_{r}$, such as ``$E_{r}^{s,t}=0$
when $s \geq 3t+4$'', then the same result holds at $E_{r+k}$ for all
$k \geq 0$.  As one might imagine, such results can be quite useful.
In a typical spectral sequence, one does not have an explicit
description of $E_{r}$ when $r \geq 3$, so most vanishing theorems
have been proven at the $E_{1}$-- or $E_{2}$--term; the main result of
\cite{miller-wilkerson} is a good example of this.  One might ask how
to produce vanishing lines at later terms of a spectral sequence.  Our
main theorem provides a way of doing this; see the second and third
applications after the statement of the theorem for examples.

The following definition is at the heart of the paper.

\begin{definition}
A property $P$ of spectra is \emph{generic} if
\begin{itemize}
\item whenever a spectrum $X$ satisfies $P$, then so does
any retract of $X$; and
\item if $X \rightarrow Y \rightarrow Z$ is a cofibration and two of
$X$, $Y$ and $Z$ satisfy $P$, then so does the third.
\end{itemize}
In other words, a property is generic if the full subcategory of all
spectra satisfying it is thick.
\end{definition}

For example, for any spectrum $E$, the property that $E_{*}(X) = 0$ is
a generic property of $X$.

We recall three other concepts: first, given a connective spectrum $W$,
we write $\conn{W}$ for its connectivity; that is, $\conn{W}$ is the
largest number $w$ so that $\pi_{n}W=0$ when $n\leq w$.  Second, we
assume that our ring spectrum $E$ satisfies the standard assumptions
for the construction and convergence of the $E$--based Adams spectral
sequence---in other words, the assumptions necessary for Theorem
15.1(iii) in \cite[Part III]{adams-blue}; see also
Assumptions~2.2.5(a)--(c) and (e) in \cite{ravenel-green}.  More
precisely, we assume that $E$ satisfies the following.

\begin{condition}\label{cond-ring}
\begin{itemize}\item[]
\item $E$ is a ring spectrum, associative and commutative up to homotopy.
\item $\pi_{r}E=0$ for $r<0$.
\item The map $\pi_{0} E \otimes \pi_{0} E \xrightarrow{} \pi_{0} E$,
induced by multiplication on $E$, is an isomorphism.
\item Let $R$ be the largest subring of the rationals to which the
ring map $\Z \xrightarrow{} \pi_{0} E$ extends; then $H_{n}(E;R)$ is
finitely generated over $R$ for all $n$.
\end{itemize}
\end{condition}

Third, in order to state the main theorem precisely, we need to work
with $E$--complete spectra: ``$E$--complete'' is a technical condition
on a spectrum $X$ which ensures that the $E$--based Adams spectral
sequence converges to $\pi_{*} X$.  We give a careful definition below
(Definition~\ref{defn-complete}); for now, we point out that if the
spectrum $X$ is connective, then $E$--completeness is a mild
restriction.  For example, if $\pi_{0}E = \Z$, then every connective
$X$ is $E$--complete; if $\pi_{0}E = \Z_{(p)}$, then every connective
$p$--local $X$ is $E$--complete.  See \cite{adams-blue},
\cite{ravenel-green}, and \cite{bousfield-spec} for more information.

This is our main result.

\begin{theorem}\label{thm-main}
Let $E$ be a ring spectrum satisfying Condition~\ref{cond-ring}, and
consider the $E$--based Adams spectral sequence $E_{*}^{**}(X)
\Rightarrow \pi_{*}(X)$.  Fix a number $m$.  The following properties
of an $E$--complete spectrum $X$ are generic\uc
\begin{itemize}
\item [\textup{(i)}] There exist numbers $r$ and $b$ so that for all $s$ and $t$
with $s \geq m(t-s) + b$, then $E_{r}^{s,t}(X) = 0$.
\item [\textup{(ii)}] There exist numbers $r$ and $b$ so that for all finite
spectra $W$ and for all $s$ and $t$ with $s \geq m(t-s-\conn{W})+b$,
then $E_{r}^{s,t}(X \Smash W) = 0$.
\end{itemize}
\end{theorem}

\begin{remark}\label{remark-main}
\begin{itemize}\item[]
\item [(a)] One usually draws Adams spectral sequences $E_{r}^{s,t}$
with $s$ on the vertical axis and $t-s$ on the horizontal; in terms of
these coordinates, the properties say that $E_{r}^{s,t}$ is zero above
a line of slope $m$, with $s$--intercept $b$ in (i), and $s$--intercept
$b-m\conn{W}$ in (ii).
\item [(b)] Assuming that $X$ is $E$--complete ensures that the
spectral sequence converges, which we need to prove the theorem.  We
do not need to identify the $E_{2}$--term of the spectral sequence, so
we do not need to know that $E$ is a flat ring spectrum---one of the
other standard assumptions on ring spectra in discussions of the Adams
spectral sequence.
\end{itemize}
\end{remark}

We mention three applications of the theorem.  Since for each prime $p$
there is a classification of the thick subcategories of the category
of finite $p$--local spectra (see \cite[Theorem 7]{hopkins-smith}),
then in this setting one may be able to identify all spectra with
vanishing line of a given slope.  Consider, for example, the classical
mod $p$ Adams spectral sequence.  This is based on the ring spectrum
$E=\HFp$, and every finite $p$--local torsion spectrum is
$\HFp$--complete.  When $p=2$, for instance, since the mod $2$ Moore
spectrum has a vanishing line of slope $\frac{1}{2}$ at the
$E_{2}$--term, then the mod $2^{n}$ Moore spectrum, and indeed any type
1 spectrum, has a vanishing line of slope $\frac{1}{2}$ at some
$E_{r}$--term.  Similarly, part of the proof of the thick subcategory
theorem is the construction of a type $n$ spectrum with a vanishing
line of slope $\frac{1}{2^{n+1}-2}$ at the $E_{2}$--term, and hence any
type $n$ spectrum has a vanishing line of slope $\frac{1}{2^{n+1}-2}$
at some $E_{r}$--term.  At odd primes, any type $n$ spectrum has a
vanishing line of slope $\frac{1}{2p^{n}-2}$ at some $E_{r}$--term of
the classical mod $p$ Adams spectral sequence.  Theorem~\ref{thm-main}
gives no control over the term $r$ or the intercept $b$ of the
vanishing line.

Since the proof of Theorem~\ref{thm-main} is formal, it also applies
in any category which satisfies the axioms of a stable homotopy
category, as given in \cite{axiomatic}.  The second author has used
this result in an appropriate category of modules over the Steenrod
algebra to prove a version of Quillen stratification for the
cohomology $\cohomology$ of the mod 2 Steenrod algebra $A$ (see
\cite{palmieri-f-iso}).  A key part of the argument is to show that
for a particular ring $R$, there is a map $\phi \mathcolon \cohomology 
\xrightarrow{} R$ satisfying these two properties:
\begin{itemize}
\item Every element in $\ker \phi$ is nilpotent.
\item For every element $y$ in $R$, there is an $n$ so that $y^{2^{n}} 
\in \im \phi$.
\end{itemize}
To prove this, one uses a certain generalized Adams spectral sequence
in the setting of $A$--modules.  This spectral sequence converges to
$\cohomology$, and the ring $R$ is the zero line of the $E_{2}$--term.
The map $\phi$ is then the edge homomorphism.  To prove that $\phi$
has the desired properties, one uses Theorem~\ref{thm-main} together
with some specific computations to show that for any $m>0$, the
spectral sequence eventually has a vanishing line of slope $m$.  Every
element in $\ker \phi$ is represented by an element above the zero
line, and hence must be nilpotent---its powers lie along a line of
some slope, and some $E_{r}$--term has a vanishing line of smaller
slope.  For every $y \in R$, one can show that for all $n$, the
targets of the possible differentials on $y^{2^{n}}$ all lie above a
certain line of positive slope, and since some $E_{r}$--term has a
vanishing line of smaller slope, then $y^{2^{n}}$ is a permanent cycle
for all sufficiently large $n$.  See \cite{palmieri-f-iso,
palmieri-steenrod} for details.

A third application is the computation of localized homotopy groups.
Suppose that $E$ is a ring spectrum satisfying
Condition~\ref{cond-ring}, and fix $v \in \pi_{*}(E)$.  Say that an
$E$--complete spectrum $X$ has a ``nice $v$--action'' if $X$ has a
self-map which induces multiplication by some power of $v$ on
$E_{*}X$, and if there is an $n$ so that the $E$--based Adams spectral
sequence converging to $\pi_{*}X$ is a spectral sequence of
$\pi_{0}(E) [v^{n}]$--modules, compatibly with the $v$--self map on $X$.
Note that multiplication by $v^{n}$ does not increase Adams
filtration: it acts ``horizontally'' in this spectral sequence.  Given
a spectrum $X$ with a nice $v$--action, if one inverts $v^{n}$ in the
$E$--based Adams spectral sequence $E_{*}^{**}(X)$, one gets a
localized spectral sequence $v^{-1} E_{*}^{**}(X)$.  A natural
question is, can one compute $v^{-1} \pi_{*} X$ from the localized
spectral sequence?  It is not hard to show that if the unlocalized
spectral sequence $E_{*}^{**}(X)$ has a horizontal vanishing line at
some $E_{r}$--term, then the localized spectral sequence will converge
to $v^{-1} \pi_{*} X$: the vanishing line forces both the
differentials and the extensions to behave well.  Hence if one can
find a spectrum $X$ for which $E_{r}^{**}(X)$ has a horizontal
vanishing line, one can conclude from Theorem~\ref{thm-main} that for
all spectra with nice $v$--actions in the thick subcategory generated
by $X$, the localized spectral sequence converges and computes the
localized homotopy groups.

\section{Composites of maps in towers}\label{sec-towers}

In this section, we prove a lemma about genericity and composites of
maps in towers.  In the next section, we will apply this lemma to the
$E$--based Adams tower to prove Theorem~\ref{thm-main}.

Let $F_{\bullet}$ be an exact functor from spectra to towers of
spectra; in other words, given a spectrum $X$, we have a diagram,
functorial in $X$:
\[
\cdots \xrightarrow{} F_{n+1} X \xrightarrow{} F_{n} X \xrightarrow{}
F_{n-1}X \xrightarrow{} \cdots.
\]
By ``exact'', we mean that if $X \xrightarrow{} Y \xrightarrow{} Z$ is
a cofibration, then so is $F_{n}X \xrightarrow{} F_{n}Y \xrightarrow{}
F_{n}Z$ for each $n$.

We want to show that the following property is generic in $X$: all
$r$--fold composites of maps in the tower $F_{\bullet} X$ are zero on
homotopy, ie, are null after composition with any map from a sphere.
Actually, we refine this condition in two ways: first, rather than
using spheres as test spectra, we use an arbitrary set of spectra; and
second, we only look at $r$--fold composites $F_{s+r}X \xrightarrow{}
F_{s}X$ when $s$ satisfies an inequality which may depend on the
particular test spectra used.

\begin{lemma}\label{lemma-towers}
Fix a number $m$ and an exact functor $F_{\bullet}$ from spectra to
towers of spectra.  Let $\test = \{(W_{\alpha}, n_{\alpha})\}$ be
a set of pairs of the form \ulp{}spectrum, number\urp{}.  The
following condition is generic in $X$\uc
\begin{itemize}
\item [$(*)$] There exist numbers $r$ and $b$ so that for all
$(W_{\alpha},n_{\alpha}) \in \test$, all $s$ with $s \geq mn_{\alpha}
+ b$, and all degree zero maps $W_{\alpha} \xrightarrow{} F_{s+r}X$,
the composite $W_{\alpha} \xrightarrow{} F_{s+r}X \xrightarrow{}
F_{s}X$ is null.
\end{itemize}
\end{lemma}

In the next section, we will apply this in the two cases 
\begin{gather*}
\test = \{(S^{n},n) \ \colon \ n \in \Z \}, \\
\test = \{(W,-\conn{DW}) \ \colon \ W \ \text{finite} \},
\end{gather*}
where $DW$ is the Spanier--Whitehead dual of $W$.

\begin{proof}
If $Y$ is a retract of $X$, then the tower $F_{\bullet} Y$
is a retract of $F_{\bullet} X$, so if $F_{s+r}X
\xrightarrow{} F_{s}X$ is null after mapping in each $W_{\alpha}$,
then so is $F_{s+r}Y \xrightarrow{} F_{s}Y$.  (Given $W_{\alpha}
\xrightarrow{} F_{s+r}Y$, then consider
\[
\begin{CD}
W_{\alpha} @>>> F_{s+r}Y @>>> F_{s}Y \\
@. @VVV  @VVV \\
@.     F_{s+r}X @>>> F_{s}X \\
@. @VVV  @VVV \\
@.  F_{s+r}Y @>>> F_{s}Y
\end{CD}
\]
Then the composite $W_{\alpha} \xrightarrow{} F_{s}Y$ factors through
$W_{\alpha} \xrightarrow{} F_{s+r}X \xrightarrow{} F_{s}X$.)

The condition $(*)$ involves numbers $r$ and $b$, and we write
$(*)_{r,b}$ if we want to specify the numbers.

Given a cofibration sequence $X \xrightarrow{f} Y \xrightarrow{g} Z$ in
which $X$ and $Z$ satisfy conditions $(*)_{r,b}$ and $(*)_{r',b'}$,
respectively, we show that $Y$ satisfies $(*)_{r+r',\max (b,
b'-r)}$.  Consider the following commutative diagram, in which the
rows are cofibrations:
\[
\begin{CD}
F_{s+r+r'} X @>>> F_{s+r+r'} Y @>{F_{s+r+r'}(g)}>> F_{s+r+r'} Z \\
   @VVV       @VV{\beta}V   @VV{\gamma}V \\
  F_{s+r} X  @>>>  F_{s+r} Y   @>>>  F_{s+r} Z \\
 @VV{\delta}V    @VV{\varepsilon}V     @VVV \\
   F_{s}X    @>>>    F_{s}Y    @>>>   F_{s}Z
\end{CD}
\]
We fix $(W_{\alpha},n_{\alpha}) \in \test$ and assume that $s \geq
mn_{\alpha} + \max (b, b'-r)$, so that we have
\begin{gather*}
s \geq mn_{\alpha} + b,\\
s + r \geq mn_{\alpha} + b'.
\end{gather*}
Given any map $\zeta \mathcolon W_{\alpha} \xrightarrow{}
F_{s+r+r'}Y$, then since $\gamma \circ F_{s+r+r'}(g) \circ
\zeta$ is null, the composite $\beta \circ \zeta$ factors through
$F_{s+r}X$, giving
\[
\zeta' \mathcolon W_{\alpha} \xrightarrow{} F_{s+r} X.
\]
Since $\delta \circ \zeta'$ is null, though, then the composite
\[
W_{\alpha} \xrightarrow{\zeta} F_{s+r+r'}Y \xrightarrow{\beta}
F_{s+r}Y \xrightarrow{\varepsilon} F_{s}Y
\]
is null.
\end{proof}

\section{Adams towers and the proof of Theorem~\ref{thm-main}}
\label{sec-proof}

The difficulty in proving a result like Theorem~\ref{thm-main} is that
the $E_{r}$--term of an Adams spectral sequence does not have nice
exactness properties if $r \geq 3$---a cofibration of spectra does not
lead to a long exact sequence of $E_{r}$--terms, for instance.  So we
prove the theorem by showing that the purported generic conditions are
equivalent to other conditions on composites of maps in the Adams
tower, and we apply Lemma~\ref{lemma-towers} to conclude that those
other conditions are generic.

We start by describing the standard construction of the Adams spectral
sequence, as found in \cite[III.15]{adams-blue},
\cite[2.2]{ravenel-green}, and any number of other places.  Given a
ring spectrum $E$, we let $\overline{E}$ denote the fiber of the unit
map $S^{0} \xrightarrow{} E$.  For any integer $s \geq 0$, we let
\begin{gather*}
F_{s} X = \overline{E}^{\Smash s} \Smash X, \\
K_{s} X = E \Smash \overline{E}^{\Smash s} \Smash X.
\end{gather*}
We use these to construct the following diagram of cofibrations, which
we call the \emph{Adams tower for $X$}:
\[
\begin{CD}
X @= F_{0} X @<{g}<< F_{1} X @<{g}<< F_{2} X @<{g}<< \cdots. \\
@. @VVV @VVV @VVV \\
@. K_{0} X @. K_{1} X @. K_{2} X
\end{CD}
\]
This construction satisfies the definition of an ``$E_{*}$--Adams
resolution'' for $X$, as given in \cite[2.2.1]{ravenel-green}---see
\cite[2.2.9]{ravenel-green}.  Note also that $F_{s} X = X \Smash F_{s}
S^{0}$, and the same holds for $K_{s} X$---the Adams tower is
functorial and exact.

We pause to define $E$--completeness.

\begin{definition}\label{defn-complete}
A spectrum $X$ is \emph{$E$--complete} if the inverse limit of its
Adams tower is contractible.
\end{definition}

Given the Adams tower for $X$, if we apply $\pi_{*}$, we get an exact
couple and hence a spectral sequence.  This is called the
\emph{$E$--based Adams spectral sequence}.  More precisely, we let
\begin{gather*}
D_{1}^{s,t} = \pi_{t-s} F_{s} X, \\
E_{1}^{s,t} = \pi_{t-s} K_{s} X.
\end{gather*}
If we let $g \mathcolon F_{s+1}X \xrightarrow{} F_{s}X$ denote the
natural map, then $g_{*} = \pi_{t-s}(g)$ is the map $D_{1}^{s+1,t+1}
\xrightarrow{} D_{1}^{s,t}$.  Then we have the following exact couple
(the pairs of numbers indicate the bidegrees of the maps):
\begin{center}
\exactcouple{D_{1}^{*,*}}{E_{1}^{*,*}}{D_{1}^{*,*}}
   {(0,0)}{(1,0)}{(-1,-1)}
\end{center}
This leads to the following $r$th derived exact couple, where
$D_{r}^{s,t}$ is the image of $g_{*}^{r-1}$, and the map
$D_{r}^{s+1,t+1} \xrightarrow{} D_{r}^{s,t}$ is the restriction of
$g_{*}$:
\begin{center}
\exactcouple{D_{r}^{*,*}}{E_{r}^{*,*}}
 {D_{r}^{*,*}} 
 {(r-1,r-1)}{(1,0)}{(-1,-1)} 
\end{center}

Unfolding the $r$th derived exact couple leads to the following exact
sequence:
\begin{equation}\label{eqn-ass-exact}
\cdots \xrightarrow{} E_{r}^{s,t+1} \xrightarrow{} D_{r}^{s+1,t+1}
\xrightarrow{} D_{r}^{s,t} \xrightarrow{} E_{r}^{s+r-1,t+r-1}
\xrightarrow{} \cdots .
\end{equation}

Fix a number $m$.  With respect to the $E$--based Adams spectral
sequence $E_{*}^{**}(-)$, we have the following conditions on a
spectrum $X$:
\begin{itemize}
\item [(1)] There exist numbers $r$ and $b$ so that for all $s$ and
$t$ with $s \geq m(t-s) + b$, then $D_{r}^{s,t}(X) = 0$.  (In other
words, the map $g_{*}^{r-1} \mathcolon \pi_{t-s}(F_{s+r-1}X)
\xrightarrow{} \pi_{t-s}(F_{s}X)$ is zero.  In other words, for all
maps $f \mathcolon S^{t-s} \xrightarrow{} F_{s+r-1}X$, the composite
$g^{r-1} \circ f$ is null.)
\item [(2)] There exist numbers $r$ and $b$ so that for all $s$ and
$t$ with $s \geq m(t-s) + b$, then $E_{r}^{s,t}(X) = 0$.
\item [(3)] There exist numbers $r$ and $b$ so that for all finite
spectra $W$, all $s$ with $s \geq -m\conn{DW} + b$, and all maps $W
\xrightarrow{} F_{s+r-1}X$, then the composite $W \xrightarrow{}
F_{s+r-1}X \xrightarrow{} F_{s}X$ is null.  (Here, $\dual{W}$ denotes
the Spanier--Whitehead dual of $W$.)
\item [(4)] There exist numbers $r$ and $b$ so that for all finite
spectra $W$ and for all $s$ and $t$ with $s \geq
m(t-s-\conn{W}) + b$, then $E_{r}^{s,t}(X \Smash W) = 0$.
\end{itemize}
As with the condition in the proof of Lemma~\ref{lemma-towers}, each
condition involves a pair of numbers $r$ and $b$, and we write
$(1)_{r,b}$ to indicate condition~(1) with the numbers specified, and
so forth.

If $m=0$, for instance, then condition~(3) says that $F_{s+r-1}X
\xrightarrow{} F_{s}X$ is a phantom map whenever $s \geq b$.  If
$m=0$, then condition~(1) says that $F_{s+r-1}X \xrightarrow{} F_{s}X$
is a ghost map (zero on homotopy) whenever $s \geq b$.

\begin{lemma}\label{lemma-main}
Fix a spectrum $X$, and fix numbers $m$, $r$, and $b$.
We have the following implications\uc
\begin{itemize}
\item [\textup{(a)}] If $r \geq -m$, then $(1)_{r,b} \Rightarrow
(2)_{r,b+r-1}$.  If $r < -m$, then $(1)_{r,b} \Rightarrow
(2)_{r,b+m-1}$.
\item [\textup{(b)}] Suppose that $X$ is $E$--complete.  If $r \geq
1-m$, then $(2)_{r,b} \Rightarrow (1)_{r,b+m}$.  If $r < 1-m$, then
$(2)_{r,b} \Rightarrow (1)_{r,b-r+1}$.
\item [\textup{(c)}] If $r \geq -m$, then $(3)_{r,b} \Rightarrow
(4)_{r,b+r-1}$.  If $r < -m$, then $(3)_{r,b} \Rightarrow
(4)_{r,b+m-1}$.
\item [\textup{(d)}] Suppose that $X$ is $E$--complete.  If $r \geq
1-m$, then $(4)_{r,b} \Rightarrow (3)_{r,b+m}$.  If $r < 1-m$, then
$(4)_{r,b} \Rightarrow (3)_{r,b-r+1}$.
\end{itemize}
\end{lemma}

(Obviously, $(3)_{r,b} \Rightarrow (1)_{r,b}$ and $(4)_{r,b}
\Rightarrow (2)_{r,b}$, but we do not need these facts.)

\begin{proof}
As above, we write $g$ for the map $F_{s+1}X \xrightarrow{} F_{s}X$
and $g_{*}$ for the map $D_{1}^{s+1,t+1} \xrightarrow{}
D_{1}^{s,t}$, so that $D_{r}^{s,t}$ is the image of 
\[
g_{*}^{r-1} \mathcolon \pi_{t-s}F_{s+r-1}X \xrightarrow{}
\pi_{t-s}F_{s}X.
\]

(a)\qua Assume that if $s \geq m(t-s) + b$, then $D_{r}^{s,t} = 0$.  In
the case $r \geq -m$, if $s \geq m(t-s) + b$, then $s+r \geq
m((t+r-1)-(s+r)) + b$; so we see that $D_{r}^{s+r,t+r-1}=0$.  By the
long exact sequence \eqref{eqn-ass-exact}, we conclude that
$E_{r}^{s+r-1,t+r-1} = 0$ when $s \geq m(t-s)+b$.  Reindexing, we find
that $E_{r}^{p,q} = 0$ when $p \geq m(q-p) + b+r-1$; ie, condition
$(2)_{r,b+r-1}$ holds.  The case $r < -m$ is similar; in this case,
the long exact sequence implies that $E_{r}^{s-1,t} = 0$.

(b)\qua Assume that $r \geq 1-m$.  If $E_{r}^{s,t}(X) = 0$ whenever $s
\geq m(t-s) + b$, then $E_{r}^{s+r-1,t+r-2}(X) = 0$ when $s \geq
m(t-s) + b$.  So by the exact sequence \eqref{eqn-ass-exact}, we see
that $D_{r}^{s+1,t} \xrightarrow{} D_{r}^{s,t-1}$ is an isomorphism
under the same condition.  This map is induced by $g_{*} \mathcolon
\pi_{t-s-1}F_{s+1}X \xrightarrow{} \pi_{t-s-1} F_{s} X$, so we
conclude that when $s \geq m(t-s)+b$, we have
\begin{gather*}
\invlim_{q} \pi_{t-s-1} F_{q}X = D_{r}^{s,t-1}, \\
\invlim_{q}^{1} \pi_{t-s-1} F_{q}X = 0.
\end{gather*}
But since $X$ is $E$--complete, then $\invlim_{q}
\pi_{t-s-1} F_{q}X = 0$, so $D_{r}^{s,t-1} = \im g_{*}^{r-1} = 0$.
Reindexing gives $D_{r}^{p,q} = 0$ when $p \geq m(q+1-p)+b$; ie,
$(2)_{r,b}$ implies $(1)_{r,b+m}$.

If $r < 1-m$, then a similar argument shows that $D_{r}^{s-r+1,t-r+1}
= 0$.

Parts (c) and (d) are similar.  
\end{proof}

\begin{proof}[Proof of Theorem~\ref{thm-main}]
This follows immediately from Lemmas~\ref{lemma-main} and
\ref{lemma-towers}.  More precisely, to show that condition~(i) is
generic, one applies Lemma~\ref{lemma-towers} to the set $\test =
\{(S^{n},n) \ \colon \ n \in \Z \}$.  For condition~(ii), one applies
it to the set $\test = \{(W,-\conn{DW}) \ \colon \ W \ \text{finite}
\}$.
\end{proof}
\np


\begin{thebibliography}



\bibitem{adams-blue} {\bf J\,F Adams}, \emph{Stable homotopy and
generalised homology}, Chicago Lectures in Mathematics, University of
Chicago Press, Chicago--London (1974)

\bibitem{bousfield-spec} {\bf A\,K Bousfield}, \emph{The localization
of spectra with respect to homology}, Topology, 18 (1979) 257--281

\bibitem{axiomatic} {\bf M Hovey}, {\bf J\,H Palmieri}, and {\bf N\,P
Strickland}, \emph{Axiomatic stable homotopy theory}, Mem.
Amer. Math. Soc. {128} (1997), no.~610 

\bibitem{hopkins-smith} {\bf M\,J Hopkins}, {\bf J\,H Smith},
\emph{Nilpotence and stable homotopy theory {I}{I}}, Annals of Math.
{148} (1998) 1--49

\bibitem{miller-wilkerson} {\bf H\,R Miller}, {\bf C Wilkerson},
\emph{Vanishing lines for modules over the {S}teenrod algebra}, J. Pure
Appl. Algebra, {22} (1981) 293--307

\bibitem{palmieri-f-iso} {\bf J\,H Palmieri}, \emph{Quillen
stratification for the {S}teenrod algebra}, Annals of Math.
(1999) to appear

\bibitem{palmieri-steenrod} {\bf J\,H Palmieri}, \emph{Stable homotopy
over the {S}teenrod algebra}, Mem. Amer. Math. Soc. to appear

\bibitem{ravenel-green} {\bf D\,C Ravenel}, \emph{Complex cobordism
and stable homotopy groups of spheres}, Academic Press, Orlando,
Florida (1986)

\end{thebibliography}
\end{document}